\documentclass[icnaam,fleqn]{w-art}
\usepackage{times}
\usepackage{w-thm}
\usepackage[]{graphicx}
\usepackage{amsmath}
\usepackage{amsfonts}

\newcommand{\field}[1]{\mathbb{#1}}

\newcommand{\R}{\field{R}}  

\newcommand{\G}{\mathcal{G}}

\begin{document}
\Year{2004}
\pagespan{3}{}
\Receiveddate{11 July 2005}
\Reviseddate{}
\Accepteddate{15 July 2005}


\title[Uncertainty and Clifford FT]{Uncertainty Principle 
       for the Clifford Geometric Algebra $Cl_{3,0}$
       based on Clifford Fourier Transform}


\author[E.M.S. Hitzer]{Eckhard S.M. Hitzer\footnote{Corresponding
     author: e-mail: {\sf hitzer@mech.fukui-u.ac.jp}, Phone: +81\,776\,27\,8514,
     Fax: +81\,776\,27\,8494}} \address[]{Dep. of App. Phys., 
     Univ. of Fukui, 3-9-1 Bunkyo, 910-8507 Fukui, Japan}
\author[B. Mawardi]{Bahri Mawardi\footnote{e-mail: 
        {\sf mawardi@quantum.apphy.fukui-u.ac.jp}}}
\begin{abstract}
\end{abstract}
\maketitle                   






\section{Introduction}

In the field of applied mathematics the Fourier transform has 
developed into an important tool. It is a powerful method for solving partial differential 
equations. The Fourier transform provides also a technique for signal analysis where
the signal from the original domain is transformed  to the spectral or frequency domain. In
the frequency domain many characteristics of the signal are revealed. With these
facts in mind, we extend the Fourier transform in geometric algebra.

Brackx et al.~\cite{BDS:CA} extended the Fourier transform to multivector valued 
function-distributions in $Cl_{0,n}$ with compact support. They also showed some properties
of this generalized Fourier transform. 
A related applied approach for hypercomplex Clifford Fourier Transformations in $Cl_{0,n}$ 
was followed by B\"{u}low et. al.~\cite{GS:GCwCA}.
In~\cite{LMQ:CAFT94}, Li et. al. extended the Fourier
Transform holomorphically to a function of \textit{m} complex variables.

In this paper we adopt and expand the generalization of the Fourier transform 
in Clifford geometric algebra\footnote{For further details and proofs 
compare~\cite{HM:CFTUP}.
In the geometric algebra literature~\cite{HS:CAtoGC}
instead of the mathematical notation $Cl_{p,q}$ the notation $\G_{p,q}$ is
widely in use. It is convention to abbreviate $\G_{n,0}$ to $\G_{n}$.}
$\G_3$ recently suggested by Ebling and Scheuermann~\cite{ES:CFTonVF}. 
We explicitly show detailed properties of the real\footnote{The meaning of
\textit{real} in this context is, that we use the three dimensional volume element
$i_3 = \mbox{\boldmath $e$}_{123}$ of the geometric algebra $\G_3$ over the field
of the reals $\R$ to construct the kernel of the Clifford Fourier transformation
of definition~\ref{df:CFT}. This $i_3$ has a clear geometric interpretation.
}
Clifford geometric algebra Fourier transform (CFT),
which we subsequently use to define and prove the uncertainty principle for 
$\G_3$ multivector functions.

\section{Clifford's Geometric Algebra $\G_3$ of $\R^3$}

Let us consider an orthonormal vector basis
$\{\mbox{\boldmath $e$}_1,\mbox{\boldmath $e$}_2,\mbox{\boldmath $e$}_3\}$
of the real 3D Euclidean vector space $\R^3$. 
The geometric algebra over $\R^{3}$ denoted by $\G_{3}$ 
then has the graded $2^3=8$-dimensional basis
\begin{equation}
\label{eqp1}
\{1,\mbox{\boldmath $e$}_1,\mbox{\boldmath $e$}_2,\mbox{\boldmath $e$}_3
, \mbox{\boldmath $e$}_{12},\mbox{\boldmath $e$}_{31}, \mbox{\boldmath $e$}_{23},
\mbox{\boldmath $e$}_{123}\}.
\end{equation}
The \textit{grade selector} is defined as $\langle M \rangle_k$ for the 
$k$-vector part of $M$, especially $\langle M \rangle=\langle M \rangle_0$.
Then $M$ can be expressed as
\begin{equation}\label{eqmd1}
M = \;\langle M \rangle + \langle M \rangle_1 + \langle M \rangle_2 + \langle M \rangle_3.
\end{equation}
The \textit{reverse} of $M$ is defined by the anti-automorphism 
\begin{equation}\label{eqmd1a}
\widetilde{M}=\; \langle M \rangle 
          + \langle M \rangle_1 
          - \langle M \rangle_2 
          - \langle M \rangle_3.
\end{equation}
The \textit{square norm} of $M$ is defined by
\begin{equation}\label{eqmd3}
\|M\|^2= \;\langle M\widetilde{M}\rangle, 
\end{equation}
where 
\begin{equation}
\label{eq:SPcf}
\langle M\widetilde{N} \rangle= M*\widetilde{N}= \sum_{A}\alpha_{A} \beta_{A} 
\end{equation}
is a real valued (inner) \textit{scalar product} for any $M, N$ in $\G_3$ 
with $M=\sum_{A}\alpha_{A} \mbox{\boldmath $e$}_{A}$ 
and  $N =\sum_{A}\beta_{A} \mbox{\boldmath $e$}_{A}$, 
$A \in \{0,1,2,3,12,31,23,123\}$, $\,\,\,\alpha_{A},\beta_{A} \in \R $.
As a consequence we obtain the 
\textit{multivector Cauchy-Schwarz inequality}
\begin{equation}\label{eq28}
|\langle M \widetilde{N} \rangle|^2
\;\leq\; \|M\|^2\;\|N\|^2\qquad \textrm{for \;all}\quad M, N \in\G_3.
\end{equation}

\section{Multivector Functions, Vector Differential and Vector Derivative}

Let $f = f(\mbox{\boldmath $x $})$ be a multivector-valued function
of a vector variable $\mbox{\boldmath $x$}$ in $\G_3$. 
For an arbitrary vector $\mbox{\boldmath $a$}$
we define\footnote{
Bracket convention: 
$A\cdot B C = (A\cdot B) C \neq A\cdot (BC)$ and
$A\wedge B C = (A\wedge B) C \neq A\wedge (BC)$
for multivectors $A,B,C \in \G_{p,q}$. 
The vector variable index \boldmath $x$ of the
vector derivative is dropped: 
$\nabla_{\mbox{\boldmath ${x}$}}=\nabla$ and
$\mbox{\boldmath $a$}\cdot\nabla_{\mbox{\boldmath $x$}}=\mbox{\boldmath $a$}\cdot\nabla $,
but not when differentiating with respect to a different vector variable 
(compare e.g. proposition~\ref{prp:dfd}).
}
the \textit{vector differential} in the $\mbox{\boldmath $a$}$ direction as
\begin{equation}
\label{eqma1}
\mbox{\boldmath $a$}\cdot\nabla f(\mbox{\boldmath $x$})=\lim_{\epsilon  \rightarrow 0}
\frac{f(\mbox{\boldmath $x$}+\epsilon \mbox{\boldmath $a$})-f(\mbox{\boldmath $x$})}{\epsilon }
\end{equation}
provided this limit exists and is well defined. 
The basis independent \textit{vector derivative} $\nabla$
defined
in~\cite{HS:CAtoGC,EH:VDC} obeys equation (\ref{eqma1})
for all vectors \mbox{\boldmath $a$} and can be expanded as 
\begin{equation}
\label{eq4}
\nabla = \mbox{\boldmath $e$}_k {\partial_k } = \mbox{\boldmath $e$}_1{\partial_1}+ 
\mbox{\boldmath $e$}_2{\partial_2} + \mbox{\boldmath $e$}_3{\partial_3},
\end{equation}

\begin{proposition}
\label{prp:lin}
$\nabla (f+g)=\nabla f + \nabla g$ (linearity).
\end{proposition}
\begin{proposition}
\label{prp:prodr}
$\nabla (fg)=(\dot{\nabla}\dot{f} )g + \dot{\nabla}f \dot{g}
=(\dot{\nabla}\dot{ f})g + \sum_{k=1}^n \mbox{\boldmath $e$}_k f ({\partial_k g})$.
\\
(Multivector functions $f$ and $g$ do not necessarily commute.)
\end{proposition}
\begin{proposition}
\label{prp:chainr}
$For f(\mbox{\boldmath $x$}) = g (\lambda (\mbox{\boldmath $x$})),
\,\lambda (\mbox{\boldmath $x$})\in \R$, 
$$\quad \mbox{\boldmath $a$} \cdot \nabla f = \{\mbox{\boldmath $a$} \cdot \nabla
\lambda(\mbox{\boldmath $x$})\} \frac{\partial g}{\partial{\lambda}}$$
\end{proposition}
\begin{proposition}
\label{prp:dfd}
$\nabla f = \nabla_{\mbox{\boldmath $\scriptstyle a$}} \,(\mbox{\boldmath $a$} \cdot \nabla f )$
\quad (derivative from differential)
\end{proposition}
Differentiating twice with the vector derivative, we get the differential Laplacian 
operator $\nabla^2$. We can write $\nabla^2 =\nabla \cdot\nabla +\nabla \wedge \nabla$.
But for integrable functions $\nabla \wedge \nabla$ = 0.  In this case we have
$\nabla^2 =\nabla \cdot\nabla.$
\begin{proposition} (integration of parts)
\label{prp:iop}
$$
\int_{\R^3} 
    g(\mbox{\boldmath$x$})
    [\mbox{\boldmath $a$}\cdot\nabla h(\mbox{\boldmath$x$})]
d^3\mbox{\boldmath$x$}
=
\left[ 
      \int_{\R^2}
          g(\mbox{\boldmath$x$})h(\mbox{\boldmath$x$})
      d^2\mbox{\boldmath $x$}
\right]^{a \cdot x =\infty}_{a\cdot x =-\infty}
\\
-\int_{\R^3}
    [\mbox{\boldmath $a$}\cdot\nabla g(\mbox{\boldmath$x$})]
    h(\mbox{\boldmath$x$})
d^3\mbox{\boldmath$x$}
$$
\end{proposition}

\section{Clifford Fourier Transform (CFT)}

\begin{defn} 
\label{df:CFT}
The Clifford Fourier transform  of $f(\mbox{\boldmath $x$})$
is the function
$\mathcal{F} \{f\}$: $\R^3 \rightarrow \G_3$ given by
\begin{equation}\label{eqmk1}
\mathcal{F}\{f\}(\mbox{\boldmath $\omega$})=\int_{\R^3} f(\mbox{\boldmath $x$})
\,e^{-i_3\mbox{\boldmath $\omega$}\cdot \mbox{\boldmath $x$}}\, d^3\mbox{\boldmath $x$},
\end{equation}
\end{defn}
where 
we can write $\mbox{\boldmath$\omega$} = \omega _1 \mbox{\boldmath $e$}_1+
\omega _2 \mbox{\boldmath $e$}_2+\omega_3\mbox{\boldmath $e$}_3$, 
$\mbox{\boldmath$x$}=x_1\mbox{\boldmath $e$}_1+x_2\mbox{\boldmath $e$}_2+
x_3\mbox{\boldmath $e$}_3$ with 
$\mbox{\boldmath $e$}_1, \mbox{\boldmath $e$}_2, \mbox{\boldmath $e$}_3$ 
the basis vectors of $\R^3$. 
Note that
\begin{equation}
\label{eq:d3x}
d^3\mbox{\boldmath $x$}= \frac{d\mbox{\boldmath $x_1$}
\wedge d\mbox{\boldmath $x_2$}\wedge d\mbox{\boldmath $x_3$}}{i_3}
\end{equation}
is scalar valued ($d\mbox{\boldmath $x_{k}$} =dx_{k}\mbox{\boldmath $e_k$}$,\;$k=1,2,3$, no summation). 
Because $i_3$ commutes with every element of $\G_3$, the Clifford
Fourier kernel $e^{-i_3\mbox{\boldmath $\omega$}\cdot \mbox{\boldmath $x$}}$
will also commute with every element of $\G_3$. 
\begin{thm}
The Clifford Fourier transform $\mathcal{F}\{f\}$ of
$f\in L^2(\R^3, \G_3)$, \\
$\int_{\R^3} \| f \|^2 d^3\mbox{\boldmath $x$} < \infty$
  is invertible and its inverse is calculated by 
\begin{equation}\label{eq11}
\mathcal{F}^{-1}[\mathcal{F}\{f\}(\mbox{\boldmath $\omega$})]= f(\mbox{\boldmath $x$})
=\frac{1}{(2\pi)^3} \int_{\R^3}\mathcal{F}\{f\}(\mbox{\boldmath $\omega$}) \, 
e^{i_3\mbox{\boldmath $\omega$}\cdot \mbox{\boldmath $x$}}\, d^3\mbox{\boldmath $\omega $}.
\end{equation}
\end{thm}
A number of properties of the CFT are listed in table \ref{tb:CFTprop}.

\begin{table}
\caption{Properties of the Clifford Fourier transform (CFT) \label{tb:CFTprop}}
\begin{center}
\begin {tabular}{llll} 
\hline
\textbf{Property}         &  \textbf{Multivector Function}      &   \textbf{CFT} 
\\ 
\hline
Linearity        & $\alpha f(\mbox{\boldmath $x$})$+$\beta \;g(\mbox{\boldmath $x$})$     &  
$\alpha \mathcal{F}\{f\}(\mbox{\boldmath $\omega$})$+ $\beta \mathcal{F}\{g\}
(\mbox{\boldmath $\omega$})$ 
\\ 
Delay    &  $f(\mbox{\boldmath $x$}-\mbox{\boldmath $a$})$   & 
$e^{-i_3\mbox{\boldmath $\omega$}\cdot  \mbox{\boldmath $a$}}
\mathcal{F}\{f\}(\mbox{\boldmath $\omega$})$    
\\ 
Shift    & $e^{i_3\mbox{\boldmath $\omega $}_0\cdot\mbox{\boldmath $x$}}f(\mbox{\boldmath $x$})$   
& $\mathcal{F}\{f\}(\mbox{\boldmath $\omega$}-\mbox{\boldmath $\omega$}_0)$ 
\\
Scaling  & $f(a\mbox{\boldmath $x$})$  &$\frac{1}{a^3}\mathcal{F}\{f\}
(\frac{\mbox{\boldmath $\omega $}}{a})$
\\
Convolution  & $(f\mbox{\boldmath $\star$}g)(\mbox{\boldmath $x$})$  &
$\mathcal{F}\{f\}(\mbox{\boldmath $\omega$})\,
\mathcal{F}\{g\}(\mbox{\boldmath $\omega $})$
\\
Vec. diff. & $\mbox{\boldmath $a$}\cdot\nabla f
(\mbox{\boldmath $x$})$  &$ i_3\;\mbox{\boldmath $a$}\cdot\mbox{\boldmath $\omega$}
\mathcal{F}\{f\}(\mbox{\boldmath $\omega$)} $
\\
&$\mbox{\boldmath $a$}\cdot\mbox{\boldmath $x$}\;
f(\mbox{\boldmath $x$})$  &$ i_3\;\mbox{\boldmath $a$}\cdot 
\nabla_{ \mbox{\boldmath $\scriptstyle \omega$}}\;
\mathcal{F}\{f\}(\mbox{\boldmath $\omega$})$ 
\\
&$ \mbox{\boldmath $x$} f(\mbox{\boldmath $x$})$ &$ i_3\; 
\nabla_{ \mbox{\boldmath $\scriptstyle \omega$}}\;
\mathcal{F}\{f\}(\mbox{\boldmath $\omega$})$ 
\\
Vec. deriv. & $ \nabla^m f (\mbox{\boldmath $x$})$  &$
(i_3\;\mbox{\boldmath $\omega$})^m\mathcal{F}\{f\}(\mbox{\boldmath $\omega$)}$
\\
Plancherel T. & 
$\langle f_1(\mbox{\boldmath $x$}) \widetilde{f_2(\mbox{\boldmath $x$})}\rangle_V=$  &
$\frac{1}{(2\pi)^3}
\langle \mathcal{F}\{f_1\}(\mbox{\boldmath $\omega$})\widetilde{\mathcal{F}
\{f_2\}(\mbox{\boldmath $\omega$})}\rangle_V$ 
\\
sc. Parseval T. & $\int_{\R^3}\,\|f(\mbox{\boldmath $x$})\|^2\,d^3\mbox{\boldmath $x$}=$  &
$\int_{\R^3}\,\|\mathcal{F}\{f\}(\mbox{\boldmath $\omega$})\|^2\,
d^3\mbox{\boldmath $\omega$}$
\end{tabular}
\end{center}
\end{table}

\section{Uncertainty Principle}

The uncertainty principle plays an important role in the development and 
understanding of quantum physics. It is also central for information processing~\cite{JR:HWI}. 
In quantum physics it states e.g. 
that particle momentum and position cannot be simultaneously known. 
In Fourier analysis such conjugate entities correspond to a
function and its Fourier transform which cannot both be simultaneously sharply localized. 
Futhermore much work 
(e.g.~\cite{JR:HWI,JC:UP}) 
has been devoted 
to extending the uncertainty principle to a function 
and its Fourier transform.  From the view point of geometric algebra an uncertainty 
principle  gives us information about how a multivector valued function 
and its Clifford Fourier transform are related.
\begin{thm}
\label{th:adotb}
Let $f$ be a multivector valued function in $\G_3$ which has
the Clifford Fourier transform $\mathcal{F}\{ f \}(\mbox{\boldmath $\omega$})$.
Assume $\int_{\R^3}\| f(\mbox{\boldmath $x$})\|^2 \;
d^3 \mbox{\boldmath $x$} = F <\infty $, then the following inequality holds 
for arbitrary constant vectors \mbox{\boldmath $a$}, \mbox{\boldmath $b$}:
\begin{equation}
\label{eqmkz}
\int_{\R^3}(\mbox{\boldmath $a$}\cdot\mbox{\boldmath $x$})^2
\|f(\mbox{\boldmath $x$})\|^2 \; d^3 \mbox{\boldmath $x$}  
\int_{\R^3}(\mbox{\boldmath $b$}\cdot\mbox{\boldmath $\omega$})^2\;
\|\mathcal{F}\{f\}(\mbox{\boldmath $\omega$})\|^2 d^3\mbox{\boldmath $\omega$}
\geq (\mbox{\boldmath $a$}\cdot\mbox{\boldmath $b$})^2\;\mbox{} \frac{(2\pi)^3}{4}F^2
\end{equation}
\end{thm}

Choosing \mbox{\boldmath $b$} = \mbox{\boldmath $\pm a$}, 
with $\mbox{\boldmath $a$}^2 = 1$ we get the following 
\textbf{uncertainty principle}, i.e.
\begin{equation}\label{eqmm}
\int_{\R^3}(\mbox{\boldmath $a$}\cdot\mbox{\boldmath $x$})^2\;
\|f(\mbox{\boldmath $x$})\|^2 \; d^3 \mbox{\boldmath $x$}  
\;\int_{\R^3}(\mbox{\boldmath $a$}\cdot\mbox{\boldmath $\omega$})^2\;
\|\mathcal{F}\{f\}(\mbox{\boldmath $\omega$})\|^2 d^3\mbox{\boldmath $\omega$}
\geq 
\frac{(2\pi)^3}{4}\; F^2.
\end{equation}
In (\ref{eqmm}) equality holds for \textit{Gaussian} multivector valued functions  
\begin{equation}\label{eqf}
f(\mbox{\boldmath $x$}) = C_0\; e^{-k\;\mbox{\boldmath $x$}^2}
\end{equation}
where $C_0 \in \G_3$ is a constant multivector, $0 < k \in \R.$
\begin{thm}
\label{th:adotb0}
For $\mbox{\boldmath $a$}\cdot\mbox{\boldmath $b$}$ =0, we get
\begin{equation}
\int_{\R^3}(\mbox{\boldmath $a$}\cdot\mbox{\boldmath $x$})^2\;
\|f(\mbox{\boldmath $x$})\|^2 \; d^3 \mbox{\boldmath $x$}  
\;\int_{\R^3}(\mbox{\boldmath $b$}\cdot\mbox{\boldmath $\omega$})^2\;
\|\mathcal{F}\{f\}(\mbox{\boldmath $\omega$})\|^2 d^3\mbox{\boldmath $\omega$}\;
\geq 0.
\end{equation}
\end{thm}

\begin{thm}
Under the same assumptions as in theorem \ref{th:adotb}, we obtain
\begin{equation}
\int_{\R^3}\mbox{\boldmath $x^2$}\;\|f(\mbox{\boldmath $x$})\|^2 \; 
d^3 \mbox{\boldmath $x$}  \;\int_{\R^3}\mbox{\boldmath $\omega^2$}\;\|\mathcal{F}\{f\}
(\mbox{\boldmath $\omega$})\|^2 d^3\mbox{\boldmath $\omega$}
\geq 3\;\frac{(2\pi)^3}{4} F^2.
\end{equation}
\end{thm}

\section{Conclusions}

The (real) Clifford Fourier transform  extends
the traditional Fourier transform on scalar functions to
multivector functions. Basic properties and rules for differentiation, convolution,
the Plancherel and Parseval theorems are demonstrated. We then presented an uncertainty 
principle in the geometric algebra $\G_3$
which describes how a multivector-valued function and its Clifford Fourier
transform relate. The formula of the uncertainty principle in $\G_3$
can be extended to $\G_n$ using properties of the Clifford Fourier transform for geometric algebras
with unit pseudoscalars squaring to -1.

It is known that the Fourier transform is successfully applied to solving physical equations
such as the heat equation, wave equations, etc. 
Therefore in the future, we can apply geometric algebra
and the Clifford Fourier transform to solve  such problems involving scalar, vector, bivector and
pseudoscalar fields.

\begin{acknowledgement}
  This research was financially supported by the Global Engineering Program
for International Students 2004 of the University of Fukui. 
We thank A. Hayashi for his continuous support 
and O. Yasukura and his colleagues from the Department of Applied Physics
of the University of Fukui
for their comments. We thank U. K\"{a}hler for his hints 
and E.J. Bayro-Corrochano and H. Ishi for good discussions. 
E. Hitzer wants to thank his wife and children for their patient support and
W. Sproessig for the invitation to ICNAAM 2005.
Soli Deo Gloria.
\end{acknowledgement}

\end{document}